\documentclass[11pt,reqno]{amsart}

\usepackage{latexsym,color,amsmath,amsthm,amssymb,amscd,amsfonts,graphicx,xcolor}

\newtheorem{theorem}{Theorem}
\newtheorem{corollary}{Corollary}

\newtheorem{lemma}{Lemma}

\newtheorem{proposition}{Proposition}

\def\R{{\mathbb R}}

\def\conv{{\rm conv}}

\def\span{{\rm span}}
\def\dim{{\rm dim}}

\begin{document}

\title{On the volume of sections of a convex body by cones}

\author{Matthieu Fradelizi, Mathieu Meyer, and Vlad Yaskin}

\thanks{The third author is supported in part by NSERC.  Part of this work was done when the third author was visiting Universit\'e Paris-Est Marne-la-Vall\'ee. He is grateful for its hospitality.}

\address{Universit\'e Paris-Est, Laboratoire d'Analyse et de Math\'ematiques Appliqu\'ees UMR 8050, UPEMLV, UPEC, CNRS
F-77454, Marne-la-Vall\'ee, France } \email{matthieu.fradelizi@u-pem.fr}

\address{Universit\'e Paris-Est, Laboratoire d'Analyse et de Math\'ematiques Appliqu\'ees UMR 8050, UPEMLV, UPEC, CNRS
F-77454, Marne-la-Vall\'ee, France } \email{mathieu.meyer@u-pem.fr}

\address{Department of Mathematical and Statistical Sciences, University of Alberta, Edmonton, Alberta T6G 2G1, Canada} \email{yaskin@ualberta.ca}

\date{\today}
\subjclass[2010]{Primary 52A20, 52A40}

\keywords{Convex body, section, $s$-concave function, centroid}

\maketitle

\begin{abstract} Let $K$ be a convex body in $\R^n$. We prove that in small codimensions, the sections of a convex body through the centroid
are quite symmetric with respect to volume. As a consequence of our estimates
we give a positive answer to a problem posed by M.~Meyer and S.~Reisner regarding convex intersection bodies.
\end{abstract}

\section{Introduction}\label{sec:intro}
Let $K$ be a convex body in $\mathbb R^n$. The centroid of $K$ is the point  $\frac{1}{|K|} \int_K x\, dx,$ where $|K|$ denotes the volume of $K$.  According to a result of Gr\"unbaum, if the centroid of a convex body $K$ is at the origin, then for every $u\in S^{n-1}$ we have
$$|K\cap u^+| \ge \left(1+\frac{1}{n}\right)^{-n}|K|\ge e^{-1} |K|,$$
where $u^+ = \{ x\in \mathbb R^n: \langle u, x\rangle \ge 0\}$, with equality on the left hand side inequality if and only if $K$ is a cone in the direction $u$, {\em i.e.} $K=\conv(a,K\cap (u^\bot+b))$, for some $a,b\in\R^n$.

It is a natural question whether there is a similar result for sections of $K$; in other words, does there exist an absolute constant $c>0$ such that $$|K\cap v^\perp \cap u^+ |_{n-1} \ge c | K\cap v^\perp|_{n-1},$$
for every $u,v\in S^{n-1}$ such that $u\ne \pm v$? Here and below, $v^\perp = \{ x\in \mathbb R^n: \langle v, x\rangle = 0\}.$ The centroid of $K\cap v^\perp$ may not be at the origin and so we cannot apply Gr\"unbaum's result. Nevertheless, we show in this paper that the answer to the latter question is affirmative. More generally, there is an absolute constant $c>0$ such that for every  convex body $K\subset \mathbb R^n$, every $(n-k)$-dimensional subspace $V$, $0\le k\le n-1$ and any $u\in S^{n-1}\cap V$ we have
 $$ |K\cap V \cap u^+|_{n-k}  \ge \frac{c}{(k+1)^2} \left(1+\frac{k+1}{n-k}\right)^{-(n-k-2)}{|K\cap V |_{n-k}}.$$
 When working with volumes of sections, we often use a subindex to emphasize the dimension of the corresponding subspace, even though most of the time the dimension  will be clear from the context.

 The results presented in this paper are actually proved in a more general setting, than described above. More specifically, we use cones instead of half-spaces; see Theorem \ref{main} below for a precise statement.

 Our results also allow us to give a positive answer to a problem posed in \cite{MR} by  Meyer and  Reisner. To state this problem, let us recall a few definitions. A {\it star body} $K$ in $\mathbb{R}^n$ is a compact set such that $[0,x]\subset K$ for every $x\in K$, and whose {\it radial function} defined by
$$
r_K(u) = \max\{a \ge 0 :a u \in K\}, \quad u\in S^{n-1},
$$
is positive and continuous. Geometrically, $r_K(u)$ is the distance from the origin to the point on the boundary in the direction of $u$. The following concept was introduced by Lutwak \cite{L}.
  The {\it intersection body of $K$}  is a star body $I(K)$ whose radial function is defined by
 $$r_{I(K)}(u) = |K\cap u^\perp|,$$
 for all $u\in S^{n-1}$.

  Busemann's theorem asserts that if $K$ is origin-symmetric and convex, then $I(K)$ is also convex. Without the symmetry assumption this statement is generally not true. In order to rectify the situation, Meyer and Reisner \cite{MR} suggested a new construction, which  allowed to extend Busemann's theorem to non-symmetric bodies. Let $K$ be a convex body whose centroid is at the origin. Define   the {\it convex interesection body
 $CI(K)$ of $K$} by its radial function:
$$r_{CI(K)}(u)= \min_{z\in u^{\perp} }|(P_u K^*)^{*z}|, \quad u\in S^{n-1},$$
where $K^*$ is the polar body of $K$ with respect to $0$, $P_u$ is the orthogonal projection onto the hyperplane $u^{\perp}$, and for
$C\subset u^{\perp}$ and $z\in u^{\perp}$, $C^{*z}$ is the polar body of $C$ in $u^{\perp}$ with respect of $z$:
$$C^{*z}=\{ y\in u^{\perp}: \langle y-z,x-z\rangle \le 1 \hbox{ for all } x\in C\}.$$
It was proved in \cite{MR} that $CI(K)$ is always convex, whenever $K$ is convex. Moreover,   $CI(K)\subset I(K)$, with equality if and only if $K$ is origin-symmetric.

Recall that $K$ is said
to be in isotropic position if the centroid of $K$ is at the origin    and the integral $\int_K \langle x, u\rangle^2\, dx$ is independent of $u \in S^{n-1}$. Note that every convex body can be put in isotropic position. It is a well-known result (see Hensley \cite{H}, Ball \cite{B88}, Sch\"utt \cite{Schu}, Fradelizi \cite{F99})  that for a convex body $K$ in isotropic position its intersection body $I(K)$ is ``almost" a ball, in the sense that there is a universal constant $c$ such that $$c^{-1} \le \frac{r_{I(K)}(u_1)}{r_{I(K)}(u_2)}\le c,$$
for any $u_1$, $u_2\in S^{n-1}$.

In \cite{MR} it was asked whether the same would be true for $CI(K)$. As an application of our results we show that this question has a positive answer: there exists a universal constant $c>0$ such that
$$cI(K)\subset CI(K) \subset I(K),$$
for every convex body $K$ with centroid at the origin.

In Section $2$, we establish some preliminary lemmas involving inequalities about log-concave functions. Section 3 is devoted to the statements and proofs of the main results and in Section 4 we apply them to the convex intersection bodies.

\section{Preliminary lemmas}

 We start by stating a few classical lemmas that will be needed for the proof of the main theorem.
 Recall that for every $x>0$, $$\Gamma(x)=\int_0^\infty t^{x-1}e^{-t}\,dt,$$ and
for every $x,y>0$,
$$B(x ,y) = \int_{0}^{1} t^{x-1} (1- t)^{y-1} dt=\frac{\Gamma(x)\Gamma(y)}{\Gamma(x+y)}.$$

Let $s>0$. We say that a function $f:\mathbb R^n \to \mathbb R_+$ is  $s$-concave if $f^s$ is concave on its support. Such a function is also log-concave, i.e. $\log f$ is concave. For example, if $K$ is a convex body in $\mathbb R^n$ and $F$ is a $k$-dimensional subspace, then by Brunn's theorem, the function $f(x) = |K\cap(F+x)|_{k}$, $x\in F^\bot$,  is $1/k$-concave on its support.

The following lemma is due to Berwald \cite{Be47} and Borell \cite{Bor73};
we refer to Theorem 5.12 of the survey \cite{GNT14} and Lemma 2.2.4 of the book \cite{BGVV14} for a modern presentation.

\begin{lemma}\label{berwald}
Let $m>0$ and $g:\R_+\to \R_+$ be a $1/m$-concave   integrable function. For $p>0$  denote $I_p(g)=\left(\int_{\R_+} t^{p-1} g(t)dt\right)^{\frac{1}{p}}$. Then for $0<p\le q$, one has
  \begin{equation}\label{GPVV}
    \frac{ B(p,m+1)^{\frac{1}{p}}  }{  B(q,m+1)^{\frac{1}{q}}  } g(0)^{  \frac{1}{p}-\frac{1}{q}  }
             \le \frac{ I_p(g)}{I_q(g)}\le
           \frac{  q^{\frac{1}{q}}  } {p^{\frac{1}{p}}  }\max(g)^{  \frac{1}{p}-\frac{1}{q}  } .
           \end{equation}
\end{lemma}

Recall that the Minkowski functional of a star body $K$ is defined by
$$\|x\|_K = \min\{a\ge 0: x\in aK\}.$$
If $x\in S^{n-1}$, then $\|x\|_K^{-1} = r_K(x)$.

The next lemma is mainly due to Ball \cite{B88}; see also Proposition 2.5.3, Theorem 2.5.5, and Proposition 2.5.7 of the book \cite{BGVV14}. Observe that some of our lemmas are closely related to the techniques used by O. Gu\'edon and E. Milman in \cite{GM11}.

\begin{lemma}\label{lemmaball}
Let $f:\R^k\to\R_+$ be a log-concave integrable function with $0$ in the interior of its support. For  $p> 0$ and $x\in \R^k\setminus \{0\}$, let
$$I_p(f,x)= \left(\int_{\R_+} t^{p-1} f(xt)dt\right)^{\frac{1}{p}}.$$
Then
 \vskip 0.2mm \noindent
{\bf 1)}  $x \mapsto I_p(f,x)$ is the radial function of a convex body $L_p(f)$ in $\R^k$,  with $0$ in its interior, and whose Minkowski functional is defined by
$$\|x\|_{L_p(f)}=I_p(f,x)^{-1}=\left(\int_{\R_+} t^{p-1} f(xt)dt\right)^{-\frac{1}{p}}.$$
\vskip 0.2mm \noindent
{\bf 2)} For every $u\in\R^k$, and every integer $p\ge0$ one has
$$\int_{L_{k+p}(f)}\langle x,u\rangle^p dx=\frac{1}{k+p}\int_{\R^k}\langle x,u\rangle^p f(x)dx.$$
Thus, if $$\int_{\R^k}\langle x,u\rangle f(x)dx=0, \quad \mbox{ for all } u\in S^{k-1},$$ then the centroid of $L_{k+1}(f)$ is at $0$.

Furthermore,  if for some $c>0$ and for  all $u\in S^{k-1}$ one has
$$\int_{\R^k}\langle x,u\rangle^2 f(x) dx=c,$$  then
$$\int_{L_{k+2}(f)} \langle x,u\rangle^2 dx=\frac{c}{k+2},\quad \mbox{ for all } u\in S^{k-1}.$$
\vskip 0.2mm \noindent
\end{lemma}

As a consequence of Lemmas \ref{berwald} and \ref{lemmaball}, we deduce that if a function $f$ satisfies the hypothesis of Lemma \ref{lemmaball} and is $1/m$-concave for some $m>0$, then for $0<p\le q$  one has
\begin{equation}\label{inclusions}
 \frac{ B(p,m+1)^{\frac{1}{p}}  }{  B(q,m+1)^{\frac{1}{q}}  } f(0)^{  \frac{1}{p}-\frac{1}{q}  }   L_q(f)
            \subset L_p(f)\subset
           \frac{  q^{\frac{1}{q}}  } {p^{\frac{1}{p}}  } \max(f)^{  \frac{1}{p}-\frac{1}{q}  }L_q(f).
\end{equation}

The next lemma was proved in \cite{F97}.

\vskip 3mm \noindent
\begin{lemma}\label{lemfrad}
Let $m>0$ and $f:\R^k\to\R_+$ be a $1/m$-concave integrable function such that $\int_{\R^k} \langle x,u\rangle f(x) dx=0$ for all $u\in \R^k$, then
\begin{equation}\label{centroid}
\max_{x\in \mathbb R^k} f(x) \le\left(1+\frac{k}{m+1}\right)^m f(0).\end{equation}
\end{lemma}

Let $K$ and  $L$ be convex bodies containing the origin in their interiors. The {\it geometric distance} between $K$ and $L$ is defined by
$$d_g(K,L)=\inf\{\alpha\beta: \alpha>0, \beta >0,  \frac{1}{\alpha} K\subset L\subset \beta K\}.$$
 The preceding lemmas give the following estimates.

\begin{lemma}\label{lemmaLpLk}
Let $k\ge 1$ be an integer, and $m>0, p\ge1$   with $p\le k+1$. Let $f:\R^k\to\R_+$ be a $1/m$-concave integrable function that satisfies $$\int_{\R^k} \langle x,u\rangle f(x) dx=0$$ for all $u\in \R^k$. Then one has
\begin{eqnarray}\label{geomdist}
d_g(L_{k+1}(f), L_p(f))
\le   \left(1+\frac{k}{m+1}\right)^\frac{m}{p} \frac{ ((k+1)B(k+1,m+1))^{\frac{1}{k+1}} }{\left(pB(p,m+1)\right)^{\frac{1}{p}}}.
\end{eqnarray}
Furthermore,
\begin{equation}\label{inclusion2} d     L_{k+2}(f)\subset f(0)^{-\frac{1}{(k+1)(k+2)}}L_{k+1}(f)     \subset  c e^{\frac{1}{k}}  L_{k+2}(f),\end{equation}
for some absolute constants $c,d>0$.
\end{lemma}

\noindent
{\bf Proof.}
Using (\ref{inclusions}) with $q= k+1$ and applying the bound from Lemma \ref{lemfrad}, we get
$$d_g(L_{k+1}(f), L_p(f))
\le  \left(1+\frac{k}{m+1}\right)^{\frac{m}{p}-\frac{m}{k+1}} \frac{ ((k+1)B(k+1,m+1))^{\frac{1}{k+1}} }{\left(pB(p,m+1)\right)^{\frac{1}{p}}}.
$$
To obtain (\ref{geomdist}), observe that
\begin{equation}\label{bounds}
e^{-1} \le\left(1+\frac{k}{m+1}\right)^{-\frac{m}{k+1}}\le 1.
\end{equation}
To prove (\ref{inclusion2}), we use (\ref{inclusions}) with $p = k+1$ and $q = k+2$. The conclusion follows from the left-hand side of inequality (\ref{bounds}) and the fact that there exist  universal constants $c_1,c_2>0$ such that for all $x,y\ge1$
$$\frac{c_1 x}{x+y}\le  B(x,y) ^\frac{1}{x}\le \frac{c_2 x}{x+y}.$$
Introducing new constants $c$ and $d$, we get the result.
\qed

The next lemma is well known and very old (see e.g. \cite[p. 57]{BF87}).
\begin{lemma}\label{centroidnegative}
Let $L$ be a convex body in $\R^k$ with centroid at $0$. Then $ -L\subset k L$, where  $-L=\{-x : x\in L\}$,
  and thus for all $\theta\in S^{k-1}$  one has
\begin{equation}\label{centroids}
\frac{1}{k}\le \frac{||\theta||_L}{||-\theta||_L}\le k.
\end{equation}
\end{lemma}
Another equivalent way of stating the preceding Lemma is to write that $d_g(L,-L)\le k^2$.
Using Lemmas \ref{lemmaball},  \ref{lemmaLpLk}, and \ref{centroidnegative}, we deduce the following.

\begin{lemma}\label{lemmaLpnegative}
Let $k\ge 1$ be an integer, $m>0$ and $0\le p\le k+1$. Let $f:\R^k\to\R_+$ be a $1/m$-concave integrable function that satisfies $$\int_{\R^k} \langle x,u\rangle f(x) dx=0$$ for all $u\in \R^k$. Then
\begin{equation*}
-L_p(f)\subset k\left(1+\frac{k}{m+1}\right)^\frac{m}{p} \frac{ ((k+1)B(k+1,m+1))^{\frac{1}{k+1}} }{\left(pB(p,m+1)\right)^{\frac{1}{p}}}L_p(f).
\end{equation*}
\end{lemma}

\noindent
{\bf Proof.} Note that the centroid of $L_p(f)$ is not necessarily at the origin, and so Lemma \ref{centroidnegative} does not apply. However, by  Lemma \ref{lemmaball}  the centroid of $L_{k+1}(f)$ is at the origin and, therefore,  Lemma \ref{centroidnegative} gives  $$d_g\left(-L_{k+1}(f),L_{k+1}(f)\right)\le k^2.$$
 Now we use the (multiplicative) triangle inequality  to get
\begin{eqnarray*}
& & d_g(-L_p(f),L_p(f))  \\
&&\qquad\qquad \le  d_g(-L_p(f), -L_{k+1}(f))d_g(-L_{k+1}(f),L_{k+1}(f))d_g(L_{k+1}(f),L_p(f))\\
&&\qquad\qquad \le  k^2 d_g(L_{k+1}(f),L_p(f))^2.
\end{eqnarray*}
Finally, we apply inequality (\ref{geomdist}) from  Lemma   \ref{lemmaLpLk}.
\qed

Recall that the support function of a convex body $L$ in $\R^k$ is defined by
$$h_L(u)=\sup_{x\in L}\langle x,u\rangle\quad \hbox{for $u\in\R^k$}.$$
The next lemma is a result due to \cite{KLS95}; see also \cite{FPS12} for  a simpler proof.
\begin{lemma}\label{KLS}
Let $L$ be a convex body in $\R^k$ with centroid at $0$ and let $u\in\R^k$. Then
\begin{eqnarray*}
\frac{h_L(u)^2}{k(k+2)}\leq \frac{1}{|L|}\int_L \langle x,u\rangle^2dx\leq \frac{k}{k+2}  h_L(u)^2.
\end{eqnarray*}
\end{lemma}
Notice that if we apply the preceding Lemma to $-u$ and $u$, we get that $h_L(-u)\le k h_L(u)$, which proves that $-L\subset kL$.

Using Lemma \ref{KLS}, we deduce the following proposition.

 \begin{proposition}\label{prop-almost-isot}
Let $L$ be a convex body in $\R^k$ with centroid at $0$.
 If for some  $\gamma>0,r>1$, one has $\gamma\le \int_L \langle x,u\rangle^2 dx \le \gamma r^2$ for all  $u\in S^{k-1}$, then  for $\beta(L)=\sqrt{\frac{\gamma}{|L|}}\sqrt{\frac{k+2}{k}}$  the following holds:
$$\beta(L)B_2^k\subset L\subset rk\beta(L)B_2^k.$$
The latter implies that $d_g(L, B_2^k)\le rk$, i.e. for all $\theta,\theta'\in S^{k-1}$, one has
\begin{equation}\label{isotropies}
\frac{1}{rk}\le \frac{||\theta||_L}{||\theta'||_L}\le rk\end{equation}
\end{proposition}

\begin{proposition}\label{distLkball}
Let $k\ge 1$ be an integer and let $f:\R^k\to\R_+$ be a $1/m$-concave  integrable function that satisfies $$\int_{\R^k} \langle x,u\rangle f(x) dx=0$$ for all $u\in \R^k$.
If for some  $\gamma>0,r>1$ one has $\gamma\le\int_{\R^k} \langle x,u\rangle^2 f(x) dx\le \gamma r^2$, then
\begin{equation*}
d_g(L_{k+1}(f), B_2^k)\le  r  a^k
\end{equation*}
for some universal constant $a$.
\end{proposition}

\noindent
{\bf Proof.}
By Lemma \ref{lemmaball},   the centroid of $L_{k+1}(f)$ is at $0$ and
$$ \int_{\mathbb R^k}\langle x,u\rangle^2 f(x)dx=(k+2)\int_{L_{k+2}(f)}\langle x,u\rangle^2 dx.$$
Now,  using lemma \ref{lemmaLpLk},  one has for some absolute constants $c,d>0$,
$$ d       L_{k+2}(f)\subset f(0)^{-\frac{1}{(k+1)(k+2)}}L_{k+1}(f)     \subset c e^{\frac{1}{k}}  L_{k+2}(f).$$
Thus, integrating, we get
$$d^{k+2}\int_{L_{k+2}(f)}\langle x,u\rangle^2dx\le f(0)^{-\frac{1}{k+1}}\int_{L_{k+1}(f)}\langle x,u\rangle^2dx\le c^{k+2} e^{\frac{k+2}{k}}\int_{L_{k+2}(f)}\langle x,u\rangle^2dx.$$
This means that
$$ d^{k+2}\frac{\gamma}{k+2}\le f(0)^{-\frac{1}{k+1}}\int_{L_{k+1}(f)}\langle x,u\rangle^2dx\le c^{k+2} e^{\frac{k+2}{k}}\frac{\gamma r^2}{k+2}.$$
   It follows that  $L_{k+1}(f)$ satisfies the hypotheses of Proposition \ref{prop-almost-isot} with the constants $\gamma(f)=d^{k+2}\frac{\gamma}{k+2} f(0)^\frac{1}{k+1}$ and $r(f)= e^{\frac{k+2}{2k}} (cd^{-1})^\frac{k+2}{2}r$. Therefore,
$$d_g(L_{k+1}(f), B_2^k)\le rk e^{\frac{k+2}{2k}} (cd^{-1})^\frac{k+2}{2}\le   r a^k,$$
for some universal constant $a$. \qed

\section{Main results}

\begin{theorem}\label{main}
Let $n, p, k\ge1$ be integers such that $p\le k\le n$. Let $F$ be an $(n-k)$-dimensional subspace of $\R^n$ and $C\subset F^\bot$ be a closed cone with vertex at $0$. Let $G=\span(C)$, $p=\dim(G)$, and assume that $|C\cap B_2^n|_p>0$. Let $K$ be a convex body in $\R^n$ with centroid at the origin. \\
(1) Then
 $$\frac{|K\cap(F-C)|_{n-k+p}}{|K\cap(F+C)|_{n-k+p}}
 \le k^p\left(1+\frac{k}{n+1-k}\right)^{n-k}{n+p-k\choose p} {n+1\choose k+1}^{-\frac{p}{k+1}} .$$\\
(2) If, moreover, $K$ is in isotropic position, then
$$b^{-1} \frac{|B_2^n\cap(F+C)|_{n-k+p}}{|B_2^n\cap(F+G)|_{n-k+p}}\le \frac{|K\cap(F+C)|_{n-k+p}}{|K\cap(F+G)|_{n-k+p}}\le  b \frac{|B_2^n\cap(F+C)|_{n-k+p}}{|B_2^n\cap(F+G)|_{n-k+p}},$$
where $$b = \min\left\{n^p, a^{kp}\left(1+\frac{k}{n+1-k}\right)^{n-k}{n+p-k\choose p} {n+1\choose k+1}^{-\frac{p}{k+1}}\right\},$$
and $a$ is the absolute constant from Proposition \ref{distLkball}.

\end{theorem}

\noindent
{\bf Proof.} We will first consider the case $p\le k\le n-1$. Define $f: F^\bot\to\R_+$ by  $f(x)=|K\cap(F+x)|_{n-k}$. Then the function $f$ is $1/(n-k)$-concave and
\begin{eqnarray*}|K\cap(F+C)|_{n-k+p}&=&\int_C f(x)\,dx=\int_{C\cap S^{p-1}(G)} \int_0^{+\infty}r^{p-1}f(r\theta)\,drd\theta\\
&=&\int_{C\cap S^{p-1}(G)} \|\theta\|_{L_p(f)}^{-p}\,d\theta.
\end{eqnarray*}

\noindent
{\bf Part (1).} Using  Lemma \ref{lemmaLpnegative} with $m=n-k$, we get, for every $\theta\in S^{p-1}(G)$,
$$ \|-\theta\|^{-1}_{L_p(f)}\le
k\left(1+\frac{k}{n+1-k}\right)^\frac{n-k}{p}\frac{((k+1)B(k+1,n+1-k))^{\frac{1}{k+1}} }{\left(pB(p,n+1-k)\right)^{\frac{1}{p}}} \|\theta\|^{-1}_{L_p(f)}.$$
For positive integers $p,q$   one has $pB(p,q+1)=\frac{p!q!}{(p+q)!}={p+q\choose p}^{-1}.$ Thus the preceding inequality may be written in the following form:
$$
 \|-\theta\|^{-1}_{L_p(f)}\le k\left(1+\frac{k}{n+1-k}\right)^\frac{n-k}{p}{n+p-k\choose p}^{\frac{1}{p}}{n+1\choose k+1}^{-\frac{1}{k+1}} \|\theta\|^{-1}_{L_p(f)}.
$$
The result follows by raising this inequality to the power $p$ and integrating over $C\cap S^{p-1}(G)$.

\noindent
{\bf Part (2).}
Since $K$ is isotropic,  one has for all $u\in S^{k-1} (F^{\perp})$,
$$C_K:=\int_K \langle x,u\rangle^2 dx =\int_{F^\bot}\langle x,u\rangle^2 f(x)dx.$$
Thus the hypotheses of Proposition \ref{distLkball} are satisfied with $\gamma=C_K$ and $r=1$, and so $$d_g(L_{k+1}(f)), B_2^k)\le a^k,$$ where $a$ is the constant from Proposition \ref{distLkball}.

Using this bound, the triangle inequality and Lemma \ref{lemmaLpLk}, one has
$$
d_g(L_p(f), B_2^k)\le  d_g(L_p(f), L_{k+1}(f))\, d_g(L_{k+1}(f)), B_2^k)\le   c_{n,p,k} ,
$$
where $$c_{n,p,k} = a^k \left(1+\frac{k}{n+1-k}\right)^\frac{n-k}{p}{n+p-k\choose p}^{\frac{1}{p}}{n+1\choose k+1}^{-\frac{1}{k+1}}.$$

Now define $g: F^\bot\to\R_+$ by  $g(x)=|B_2^n\cap(F+x)|_{n-k}$. Observe that $L_p(g)$ is a ball in $F^\bot$, and
\begin{equation*}|B_2^n \cap(F+C)|_{n-k+p}= \int_{C\cap S^{p-1}(G)} \|\theta\|_{L_p(g)}^{-p}d\theta.
\end{equation*}
Since $d_g(L_p(f), B_2^k) =  d_g(L_p(f), L_p(g))\le c_{n,p,k}$, we have, for all $\theta\in S^{k-1} (F^{\perp})$,
$$\frac{1}{\alpha} \|\theta\|^{-1}_{L_p(g)} \le \|\theta\|^{-1}_{L_p(f)} \le \beta \|\theta\|^{-1}_{L_p(g)},$$
where $\alpha$ and $\beta$ are positive numbers such that $\alpha \beta =c_{n,p,k}$.
Raising this inequality to the power $p$ and integrating over $C\cap S^{p-1}(G)$, we get
$$\frac{1}{\alpha^p}  |B_2^n \cap(F+C)|_{n-k+p}\le |K \cap(F+C)|_{n-k+p}\le \beta^p |B_2^n \cap(F+C)|_{n-k+p}.$$
On the other hand, integration over $ S^{p-1}(G)$ gives
$$\frac{1}{\alpha^p}  |B_2^n \cap(F+G)|_{n-k+p}\le |K \cap(F+G)|_{n-k+p}\le \beta^p |B_2^n \cap(F+G)|_{n-k+p}.$$
After dividing the last two inequalities, we get the inequality in part (2) with the bound $(c_{n,p,k})^p$.

Note that when $k$ is large (comparable to $n$), this method does not yield a good bound. Instead, one should proceed as follows. Since $K$ is in isotropic position, we have $d_g(K,B_2^n)\le n$, which follows, for example, from Proposition~\ref{prop-almost-isot}. Therefore, for all $\theta\in S^{n-1}$ we have
$$\frac{1}{\alpha}  |\theta|_2^{-1}  \le \|\theta\|^{-1}_{K} \le \beta |\theta|_2^{-1} ,$$
where $\alpha$ and $\beta$ are positive numbers such that $\alpha \beta =n$. Raising to the power $p$, integrating over $(F+C)\cap S^{n-1} $ and then over $(F+G)\cap S^{n-1}$, and dividing the corresponding inequalities, we get
$$n^{-p} \frac{|B_2^n\cap (F+ C) |_{n-k+p}}{|B_2^n\cap(F+G)|_{n-k+p}}\le \frac{|K\cap(F+C)|_{n-k+p}}{|K\cap(F+G)|_{n-k+p}}\le n^p \frac{|B_2^n\cap(F+C)|_{n-k+p}}{|B_2^n\cap(F+G)|_{n-k+p}}.$$
Let us remark that this approach does not require $k\le n-1$, and so the proof of part (2) of the theorem is complete.

To finish the proof of part (1) for the remaining case $p\le k= n$, notice that since the centroid of $K$ is at the origin, we have $$\|-\theta\|_K^{-p}\le n^p \|\theta\|_K^{-p},$$ for all $\theta\in S^{n-1}$. Integrating this inequality over $C\cap S^{p-1}(G)$, and noting that $F=\{0\}$ in this case, we get the result.

\qed

We will now discuss some corollaries of the previous theorem.
First of all, observe that there exist universal constants $c,d>0$ such that for any integers $n\ge k\ge 1$ we have
$$d\frac{n}{k}\le{n+1\choose k+1}^\frac{1}{k+1}\le c\frac{n}{k}.$$
Therefore, part (1) of Theorem \ref{main} gives
$$\frac{|K\cap(F-C)|_{n-k+p}}{|K\cap(F+C)|_{n-k+p}} \le  \left( \frac{ck^2}{n} \right) ^p\left(1+\frac{k}{n+1-k}\right)^{n-k}{n+p-k\choose p} . $$
In particular, when  $p=1$, we get the following.

\begin{corollary}\label{cor1} There is an absolute constant $c>0$ such that for any integers $n\ge k\ge 1$, any convex body $K$  in $\R^n$ whose centroid is at the origin, any  $(n-k)$-dimensional subspace $F$ of $\R^n$, and any $\theta\in S^{n-1}\cap F^\bot$, we have
 $$\frac{|K\cap(F+\R_+\theta)|_{n-k+1}}{|K\cap(F+\R_{-}\theta)|_{n-k+1}}\le  c k^2 \left(1+\frac{k}{n-k+1}\right)^{n-k-1}\le c\left(\frac{k(n+1-k)}{n+1}\right)^2e^k.$$
\end{corollary}

In particular, when $k=2$, we obtain the following result.

\begin{corollary}\label{cor2}
 There exists a constant $c>0$ such that for every $n\ge 2$, every convex body $K$ in $\R^n$,  with centroid at $0$, and every $u,v\in S^{n-1}$ such that $v\not=\pm u$, one has
$$\frac{1}{c}
\le
\frac
{\big|\{ x\in K: \langle x,u\rangle =0,\ \langle x,v\rangle \ge 0\} \big|_{n-1}} { \big|\{ x\in K: \langle x,u\rangle =0,\ \langle x,v\rangle \le 0\}\big|_{n-1}}
\le c.$$
\end{corollary}
\vskip 3mm\noindent

More generally,   the following holds.

\begin{corollary}
Let $K$ be an isotropic convex body in $\R^n$ and $1\le p\le k\le n$. Let $E$ be a $n-k+p$-dimensional subspace of $\R^n$ and
let $u_1, \dots,u_{p}$ be pairwise orthogonal vectors in $E$ then, for some universal constant $c>0$,
$$ | K\cap\bigcap_{i=1}^p\{x\in E:\langle x,u_i\rangle\ge 0 \}|_{n-k+p}\ge\max\left\{(2n)^{-p}, e^{-ckp}\right\}|K\cap E|_{n-k+p}. $$
In particular if $k=p$ and $E=\R^n$, one has
$$ | K\cap\bigcap_{i=1}^p\{x\in \R^n:\langle x,u_i\rangle\ge 0 \}|_n\ge \max\left\{(2n)^{-p}, e^{-cp^2}\right\}  {|K|_n}.$$
\end{corollary}

\noindent
{\bf Proof.} \\
Let $F=\cap_{i=1}^p\{x\in E:\langle x,u_i\rangle= 0 \}$ and $C=\{y\in E\cap F^\bot : \langle x,u_i\rangle\ge 0, 1\le i\le p\}$.
Then $F$ is $(n-k)$-dimensional and $\cap_{i=1}^p\{x\in E:\langle x,u_i\rangle\ge 0\}=F+C$. Thus we may apply Theorem \ref{main}, part (2). The estimate follows from the following inequalities.
$$\left(1+\frac{k}{n+1-k}\right)^{n-k} \le e^k,$$
$${n+p-k\choose p} \le e^p \left(\frac{n+p-k}{p}\right)^p\le e^p \left(\frac{n}{p}\right)^p,$$
$$ {n+1\choose k+1}^{-\frac{p}{k+1}} \le d^{p} \left( \frac{k}{n} \right)^p,$$
for an absolute constant $d$.

\qed

\vskip 3mm\noindent
{\bf Remarks}

\vskip 0.5mm\noindent
{\bf 1.} We will now discuss sharpness of Corollary \ref{cor1}. When $k$ is a small fixed number, then the bound in Corollary \ref{cor1} is of constant order. What happens when $k$ is not small, but close to $n$? We will see that  in this case the bound in Corollary \ref{cor1} gives the right order.  Let $\Delta_n$ be a  regular simplex in $\R^n$ with vertices $v_1,\dots,v_{n+1}$ and centroid at $0$. For $1\le l\le n-1$, let $E_l$ be the $l$-dimensional subspace of $\R^n$ generated by $v_1,\dots v_l$. Then
it is easy to see that $E_l\cap \Delta_n=\conv(v_1,\dots v_l, f_{l+1})$, where
$$f_{l+1}= \frac{\sum_{k=l+1}^{n+1} v_k}{n+1-l}=-\frac{\sum_{k=1}^l v_k}{n+1-l}.$$
Then we have
$$\frac{ |\Delta_n\cap E_l\cap f_{l+1}^+|} {|\Delta_n\cap E_l|}= \left(\frac{1/(n+1-l)}{1/(n+1-l) +1/l}\right)^l
= \left(\frac{l}{n+1}\right)^l,$$
and, therefore,
$$\frac{ |\Delta_n\cap E_l\cap f_{l+1}^-|}{ |\Delta_n\cap E_l\cap f_{l+1}^+|}= \frac{1- \big(\frac{l}{n+1}\big)^l}{\big(\frac{l}{n+1}\big)^l}
= \left(\frac{n+1}{l}\right)^l-1.$$
Taking   $l=n-k+1$ in the latter expression and comparing it to the bound in Corollary \ref{cor1}, we see that they differ by the factor  $$c \left(\frac{k (n-k+1)} {n+1}\right)^{2}.$$
The latter is essentially constant for large values of $n$ if $k$ is fixed, or $n-k$ is fixed.

\vskip 0.5mm\noindent
{\bf 2.} The result of part (1) of Theorem \ref{main}  is  sharp for $p=k=n$. If one takes the regular centered simplex $\Delta$ in $\R^n$ and a small nonempty cone $C$ around
the direction of a vertex of the simplex, then $|C\cap K|\sim n^n |(-C)\cap K|$, when the spherical measure of $C\cap S^{n-1}$ tends to $0$.

\vskip 0.5mm\noindent
{\bf 3.}
 It may be asked what is the best constant $\alpha_n>0$ such that for any  isotropic convex body $K$ in $\R^n$ and any orthonormal basis $u_1,\dots, u_n$, one has
 $$\left(\frac{|K\cap\{x:\langle x,u_i\rangle\ge0, 1\le i\le n\}|}{|K|} \right)^{1/n} \ge \frac{\alpha_n}{2}.$$
 It is clear and follows from Corollary 3 that $\alpha_n\ge n^{-1}$. Is there a better estimate? Some calculations on the regular simplex seem to indicate that $\alpha_n$ could be of the order of $n^{-\frac{1}{2}}$.
 Observe that if $K$ is centrally symmetric and in John's position (instead of being isotropic) then the bound $n^{-\frac{1}{2}}$ trivially holds because then $B_2^n\subset K\subset\sqrt{n} B_2^n$ and the constant has the right order of magnitude as shown by the example of the cube:

Let $K=B_{\infty} ^n =[-1,1]^n$. If $n=2^m$, it is well known that one can select $n$ vertices $f_1,\dots, f_n$ of some facet of $K$ such that $(f_1,\dots,f_n)$ are pairwise orthogonal.  Let $C$ be the convex cone generated by $f_1, \dots, f_n$. Then   $$ |K\cap C |= \frac{ n^{n/2} }{n!}, $$
and $\alpha_n  \le  e  n^{-\frac{1}{2}}$.

\section{Application to the convex intersection body}
Let $K\subset \mathbb R^n$ be a convex body whose centroid is at the origin. As was mentioned in the introduction, the   convex intersection body
 $CI(K)$ of $K$ is the body with the radial function
$$r_{CI(K)}(u)= \min_{z\in u^{\perp} }|(P_u K^*)^{*z}|_{n-1}, \quad u\in S^{n-1}.$$
Since $(P_u K^*)^{*}=K\cap u^{\perp}$, and since for $z\in L$,
$$|L^{*z}|_{n-1} = \int_{L^*}\frac{1}{(1-\langle z,y\rangle)^n } dy ,$$
for a convex body $L$ in $\mathbb R^{n-1}$ (see \cite{MW}),  we get
$$r_{CI(K)}(u)=\min_{z\in P_u K^*} \int_{K\cap u^{\perp}} \frac{1}{(1-\langle z,y\rangle)^n } dy.$$
Observe that for every $z\in P_u K^*$,
$$  \int_{K\cap u^{\perp} }                               \frac{1}{(1-\langle z,y\rangle)^n}  dy
\ge \int_{K\cap u^{\perp}\cap z^+ } \frac{1}{(1-\langle z,y\rangle)^n}  dy
\ge | K\cap u^{\perp}\cap z^+|. $$
By Corollary \ref{cor2}, there is an absolute constant $c>0$ such that
$$| K\cap u^{\perp}\cap z^+| > c | K\cap u^{\perp}|.$$
Therefore, for all $u\in S^{n-1}$, $$r_{CI(K)}(u)\ge c r_{I(K)}(u).$$
On the other hand,
$$r_{CI(K)}(u)\le \int_{K\cap u^{\perp}} \frac{1}{(1-\langle z,y\rangle)^n } dy\Big|_{z=0} = r_{I(K)}(u).$$
Thus, we obtain a positive answer to the conjecture raised by Meyer and Reisner in \cite{MR}.
\vskip 3mm\noindent
\begin{theorem}
There exists an absolute constant $c>0$ such that for every $n\ge 2$ and every convex body $K$ in $\R^n$ with centroid at the origin, one has
$$cI(K)\subset CI(K) \subset I(K).$$
\end{theorem}

\end{document}